\documentclass{llncs}

\usepackage{latexsym}
\usepackage{amssymb}
\usepackage{amsmath}
\usepackage{hyperref}
\usepackage{graphicx}
\usepackage{colortbl}

\usepackage[normalem]{ulem}
\usepackage{multicol}
\usepackage[backend=biber,style=numeric,sortcites,maxbibnames=9,sorting=nty,natbib,hyperref]{biblatex}

\renewbibmacro{in:}{%
  \ifentrytype{article}
    {}
    {\bibstring{in}%
     \printunit{\intitlepunct}}}

\DeclareMathOperator{\boxi}{box} 

\DeclareMathOperator{\KG}{K}

\def\Ascr{\mathcal{A}}
\def\Bscr{\mathcal{B}}
\def\Cscr{\mathcal{C}}

\def\Fscr{\mathcal{F}}

\def\Iscr{\mathcal{I}}

\def\Oscr{\mathcal{O}}

\usepackage{marginnote}
\usepackage[textwidth=2cm]{todonotes}

\newtheorem{fact}{Fact}

\title {Boxicity and Interval-Orders:}
\subtitle{Petersen and the Complements of Line Graphs} 

\author{Marco Caoduro\inst{1}\orcidID{0000-0001-9268-884X} \and Andr\'as Seb\H{o}\inst{2}\orcidID{0000-0002-1207-888X}}

\institute{Sauder School of Business, The University of British Columbia, Vancouver, Canada \email{marco.caoduro@ubc.ca} \and CNRS, Laboratoire G-SCOP, Univ.~Grenoble Alpes, Grenoble, France \email{andras.sebo@cnrs.fr}}

\begin{document}

\maketitle

\begin{abstract}
 
    The boxicity of a graph is the smallest dimension $d$ allowing a representation of it as the intersection graph of a set of $d$-dimensional axis-parallel boxes. We present a simple general approach to determining the boxicity of a graph based on studying its ``interval-order subgraphs''.
    
    The power of the method is first tested on the boxicity of some popular graphs that have resisted previous attempts:  the boxicity of the Petersen graph is $3$,  and  more generally, that of the Kneser-graphs $K(n,2)$ is $n-2$ if $n\ge 5$, confirming a conjecture of Caoduro and Lichev [Discrete Mathematics, Vol. 346, 5, 2023]. 

    Since every line graph is an induced subgraph of the complement of $K(n,2)$, the  developed tools show furthermore that line graphs have only a polynomial number of edge-maximal interval-order subgraphs. This opens the way to polynomial-time algorithms for problems that are in general $\mathcal{NP}$-hard: for the existence and optimization of interval-order  subgraphs of line-graphs, or of interval-completions of their complement.

    \keywords{Boxicity \and Interval-orders \and Interval-completion\and Kneser-graphs \and Line Graphs}
    
\end{abstract}


\newpage

\section{Introduction}
The \emph{intersection graph} of a finite family of sets $\Fscr$ is the graph $G(\Fscr)$ having vertex-set $\{v_A: A \in \Fscr\}$ and edge-set $\{v_A v_B : A,B \in \Fscr, \ A \neq B, \textrm{ and } A\cap B \neq \emptyset\}$, where multiple occurrences of a set is allowed, and different occurrences are considered as different sets.x
The \emph{boxicity} of $G=(V,E)$, denoted by $\boxi(G)$, is the minimum dimension $d$ such that $G=G(\Bscr)$, where  $\Bscr$ is a family of axis-parallel boxes in $\mathbb R^d$.
A graph $G$ has $\boxi(G) = 0$ if and only if $G$ is a complete graph,  $\boxi(G) \leq 1$ if and only if it is an {\em interval graph}, and $\boxi(G) \le k$ if it is the intersection of $k$ interval graphs \cite{1983_Cozzens}.

Complements of interval graphs will be called here {\em interval-order graphs}, referring to the natural order of disjoint intervals.   
Interval graphs, and therefore also interval-order graphs, can be recognized in linear time \cite{1976_Booth} while determining whether a graph has boxicity at most $2$ is $\mathcal{NP}$-complete~\cite{1994_Kratochvil_NP}. 
In the language of parameterized complexity, the computation of boxicity is not in the class \emph{XP} of  problems solvable in polynomial time when the parameter boxicity is bounded by a constant.

Boxicity was introduced by Roberts  \cite{1969_Roberts} in 1969, and has been a well-studied graph parameter. 
Roberts~\cite{1969_Roberts} proved that any graph $G$ on $n$ vertices has boxicity at most $\left \lfloor \dfrac{n}{2} \right \rfloor$. 
Esperet~\cite{2016_Esperet} showed that the boxicity of graphs with $m$ edges is $\mathcal{O}(\sqrt{m\log{m}})$, while Adiga, Bhowmick, and Chandran~\cite{2010_Adiga} proved that the boxicity is $\mathcal{O}(\Delta\log^2{\Delta})$ for graphs with maximum degree $\Delta$.
In~\cite{2011_Chandran}, this latter bound was improved to $\mathcal{O}(\Delta \log{\log{\Delta}})$ in the particular case of line graphs. Further relevant results have been proved by
Scheinerman~\cite{1984_Scheinerman}, Thomassen~\cite{1986_Thomassen}, Chandran and Sivadasan~\cite{2007_Chandran}, Esperet~\cite{2013_Esperet}, and Esperet, Joret~\cite{2017_Esperet}.

\medskip

We present now an important background of our results.
Let $k$ and $n$ be two positive integers such that $n \geq 2k + 1$. The \emph{Kneser-graph} $\KG(n,k)$ is the graph with vertex-set given by all subsets of $[n]:= \{1,2,\dots,n\}$ of size $k$ where two vertices are adjacent if their corresponding $k$-sets are disjoint.
Kneser-graphs stimulated deep and fruitful graph theory, ``building  bridges'' with other parts of mathematics. For instance, Lov\'asz's proof \cite{1978_Lovasz} of Kneser's conjecture~\cite{1956_Kneser} is the  source of the celebrated ``topological method.''
This method has proven to be a powerful approach for a range of challenging combinatorial problems \cite{2003_Matousek}. 

The  study of the boxicity of Kneser-graphs $\KG(n,k)$ was initiated by Mat\v{e}j Stehl\'{i}k as a question \cite{QuestionMatej}. 
Caoduro and Lichev~\cite{2023_Caoduro} established a general upper bound of $n-2$, a lower bound of $\displaystyle n - \frac{13k^2-11k+16}{2}$  for $n\geq 2k^3 - 2k^2 + 1$, and a lower bound of $n-3$ for $k=2$ that nearly matches the upper bound.
They also conjectured that $\boxi(\KG(n,2)) = n-2$ for any $n\geq 5$. 
We establish here this conjecture. Essentially new ideas are needed already for the case $n=5$; to the best of our knowledge, all declared solutions for this particular graph finish with a computer-based case-checking  (cf.~\cite[Section 6]{2023_Caoduro}).

\begin{theorem}\label{thm:box_kneser}
	The boxicity of the Kneser-graph $\KG(n,2)$ with $n\ge 5$ is $n-2$. In particular, the boxicity of the Petersen graph $\KG(5,2)$ is $3$.
\end{theorem}

The proof of this theorem starts with a well-known rephrasing of boxicity in terms of  ``interval completion'' (Section \ref{subsec:box_int}),  
allowing a graph theory perspective. This approach reveals that each interval completion of $\KG(n,2)$ is uniquely determined by the choice of  at most $5$ of its vertices. Thus, the number of interval-completions is polynomial (at most $n^5$) in general. Even though the list is nontrivial already for small Kneser graphs, there are only four essentially different interval completions of $\KG(n,2)$. 

We introduce some additional notation and terminology.
Given a graph $G$, a graph $H$ is a subgraph of $G$ if $V(H)\subseteq V(G)$  and $E(H)\subseteq E(G)$. An \emph{interval-completion} of $G$ is an interval graph containing $G$ as a subgraph, and
the \emph{line graph} of $G$, denoted by $L(G)$, is the intersection graph of the edge-set  of $G$ as a family of sets of two vertices.
The complement of a graph $G$ is denoted by $\overline G$, and the complete graph on $n$ vertices or on the set $V$ is denoted by $K_n$, and $K_V$, respectively.
Note that $\KG(n,2)=\overline{L(K_n)}$.

Leveraging the fact that the complement of every line graph $\overline{L(G)}$ is the induced subgraph of $\overline{L(K_{V(G)})}$ we get  that the complement of {\em any   line graph has only a polynomial number of interval completions:}


\begin{lemma}\label{lem:all_minimal_completion}
    Let $G=(V, E)$ be a graph
    on $n$ vertices. Then $\overline{L(G)}$ has at most $n^5$ inclusion-wise minimal interval-completions. These can be listed in  $\Oscr(n^7)$~time. 
\end{lemma}

Computing the minimum number of edges of an interval-completion for any graph is an $\mathcal{NP}$-hard problem~\cite{1979_Garey}, called {\em Interval Graph Completion (IGC)}. However, by Lemma~\ref{lem:all_minimal_completion} we have: 

\begin{theorem}\label{thm:minimum_completion} IGC is polynomial-time solvable for complements of line graphs. 
\end{theorem}

Polynomial complexity follows for detecting bounded boxicity: 

\begin{theorem}\label{thm:box_co_line}
Let $G=(V, E)$ be a graph on $n$ vertices and $k$ a positive integer. Then it can be decided in $\Oscr(n^{5k+3})$-time whether $\boxi( \overline{L(G)} ) \leq k$. 
\end{theorem}


Summarizing,  our study of Kneser graphs $\KG(n, 2)$ has uncovered fundamental structural properties of complements of line graphs. This allows us to compute in polynomial time \emph{all} minimal interval completions of complements of line graphs, that is, to list  all ``interval-disjointness subgraphs" of ``edge-disjointness graphs''. This leads to new polynomial algorithms, for instance, the weighted version of the interval completion problem for complements of line graphs.


\medskip

The paper is structured as follows.
Section~\ref{sec:prelim} introduces the first,  well-known,  classical   properties  that play a fundamental role in the study of boxicity. It also  characterizes interval-order graphs using the orderings of their vertices.
Section~\ref{sec:interval_orders} focuses on the interval-order subgraphs of  line graphs and shows that all interval-order subgraphs of the  line graph of $K_n$ can be easily described using only five vertices of $K_n$.
Section~\ref{sec:petersen} contains the proofs of the main results. It determines the boxicity of the Petersen graph, which had been considered to be  a difficult open problem. More generally, the boxicity of the Kneser-graphs $\KG(n,2)$ for any integer $n \geq 5$ is established. 
The proof is based on elementary results and  can be generalized to solve the interval-completion problem for complements of line graphs in polynomial time.
It also leads to detecting constant boxicity, placing boxicity in the class XP for complements of line graphs.
It is yet to be proved whether this problem is $\mathcal{NP}$-hard,  or can be solved in polynomial time.
The last  section (Conclusion, Section~\ref{sec:conclusion}) offers applications, possible extensions of the methods, and open questions. 

\section{Preliminaries}\label{sec:prelim}
We use standard graph theory notation (mostly following \cite{2003_Schrijver}): given a graph $G=(V, E)$ and a vertex $v \in V$, $\delta(v)$ denotes the edges incident to $v$ and $N(v)$ the vertices adjacent to $v$. The {\em degree} of a vertex $v$ is $d(v):=|\delta(v)| = |N(v)|$. For $V' \subseteq V$, $G[V']$ is the induced subgraph of $G$ with vertex-set $V'$,  and edge-set $\{uv \ : \ uv \in E \ \textrm{and} \ u,v, \in V'\}$.

Graphs are {\em simple} in this paper, that is, they do not have loops or parallel edges. 
(These are either senseless or irrelevant to the results.) 

\subsection{Defining the boxicity using interval graphs}\label{subsec:box_int}

An \emph{axis-parallel box} in $\mathbb{R}^d$ is a Cartesian product $I_1\times I_2\times \dots \times I_d$ where each $I_i$ is a closed interval in the real line.
Axis-parallel boxes $B, B'$ intersect if and only if for all $1\leq i \leq d$ the intervals $I_i$ and $I'_i$ intersect. Roberts~\cite{1969_Roberts} formalized a set of simple but important statements assuring starting tools for the study of the boxicity: the boxicity of a graph is at most $k$ if and only if it is the intersection of $k$ interval graphs; adding a vertex or even two non-adjacent vertices the boxicity increases by at most one,   so the boxicity of a graph is at most  $\displaystyle \left \lfloor \frac{n}{2} \right \rfloor$. We will extensively use the former definition of boxicity in the following more comfortable form for us (obtained by complementation using ``de Morgan's law''): 

We say that the family of edge-sets of a graph $G=(V, E)$, $\Cscr = \{C_1, .. C_k\}$ $(C_i\subseteq E)$ is a \emph{interval-order-cover}, or a \emph{$k$-interval-order-cover} of $E$ (or $G$) if $\bigcup_{i=1}^k C_i=E$, and each $(V, C_i)$ is an interval-order graph. 
The $C_i$ will mostly be called colors and an edge in $C_i$ will be said to have \emph{color} $i$. Some edges will have several colors, which may be necessary to make each color an interval-order graph. 

\begin{lemma}[Cozzens and Roberts~\cite{1983_Cozzens}, 1983]\label{lem:Coz_Rob}
	Let $G$ be a graph. Then $\boxi(G) \leq k$ if and only if $\overline G$ has  a  $k$-interval-order-cover.
\end{lemma}

\subsection{Defining interval-order graphs using orderings on their vertices}\label{subsec:ordering}

Given a family $\Iscr$ of $n$  intervals in $\mathbb{R}$, we order them in a non-decreasing order $\sigma=I_1\ldots I_n$ of  their right end-points, and consider the  interval-order graph $G=(V, E)$, and  $V=\{v_1,\ldots, v_n\}$, where $v_i$ corresponds to $I_i$ $(i=1,\ldots n)$ and two vertices are joined if and only if the corresponding intervals are disjoint (the complement of their intersection graph).

Orienting the edges of $G$ from the vertices of smaller  index towards those of larger index, we have:
\begin{equation}\label{eq:chain_property}    
	\textit{if $i > j$, then $N^+(v_i) \subseteq N^+(v_j)$,}
\end{equation}
where the \emph{out-neighborhood} of a vertex $v_l \in V(G)$ is the set $N^+(v_l) := \{v_m : v_lv_m  \in E(G), l < m\}$.

Clearly, the series of the sizes of $N^+(v_i)$ for $i=1,\ldots, n$ is (not necessarily strictly) monotone decreasing from $d(v_1)$ to $0$.
In addition, {\em in the chain of  (\ref{eq:chain_property}), only the first at most $D$ sets are not empty, where $D$ is the maximum degree of $G$.} Indeed, suppose  $N^+(v_{i}) \ne \emptyset$ for $i\ge D +1$, and let $x\in N^+(v_{i})$. Then by (\ref{eq:chain_property}), $x$ is also in the neighborhood of all previous vertices, that is,  $N_G(x)\supseteq \{v_1,\ldots, v_{D + 1}\}$, contradicting that $D$ is the maximum degree.

By necessity, one realizes that (\ref{eq:chain_property})  actually  {\em characterizes} interval-order graphs (the proof is immediate by induction).  It is explicitly stated in Olariu's paper~\cite{1991_Olariu}:

\begin{lemma}\label{lem:co_intervals} 
	Let $G$ be an undirected graph. Then $G=(V, E)$ is an interval-order graph if and only if $V$ has an ordering  $(v_1,\ldots, v_n)$ so that orienting the edges from the vertices of smaller index towards those of larger index (\ref{eq:chain_property}) holds. 
\end{lemma}

Inspired by this characterization, we can simply construct interval-order subgraphs of an arbitrary graph, and it turns out that all interval-order subgraphs are of the form  determined by Lemma~\ref{lem:co_intervals}. 

\begin{figure}[ht]
    \centering
    \includegraphics[scale=0.63]{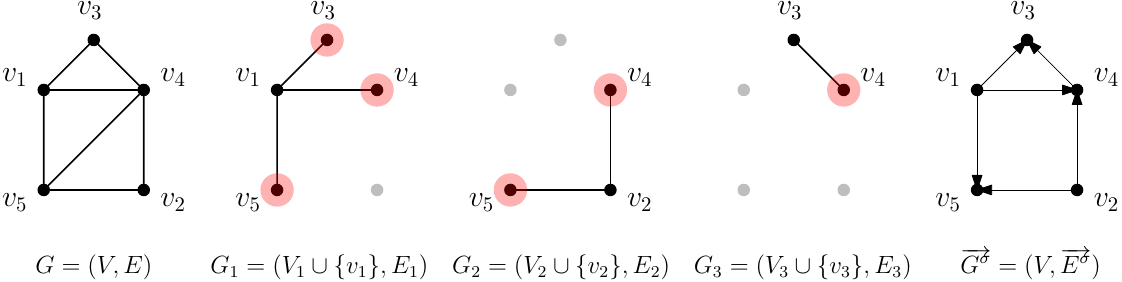}
    \caption{A graph $G$ with  an ordering $\sigma : = (v_1, v_2, v_3,  v_4, v_5)$ of $V(G)$, the corresponding $G_i$ ($i \in \{1,2,3\})$, and  $\overrightarrow{G^\sigma}$ (see Lemma~\ref{lem:co_intervals} and thereafter). For each $G_i$ ($i \in \{1,2,3\})$, the vertices of $V_i$ are marked with red disks around them. Note that $v_4$, the only vertex of $V_3$, is not contained in $N_G(v_4)$, so $V_4=\emptyset$,  $E_4 = \emptyset$.}
    \label{fig:example_G_sigma}
\end{figure}

Let $G=(V,E)$ be a graph and $\sigma$ an ordering of $V(G)$. We define the graph $G^\sigma=(V,E^\sigma)$ as follows: let $V_0:=V$, $V_i := V_{i-1}\cap N_G (v_i)$ for $1\leq i \leq n$, and $E^\sigma:= E_1 \cup E_2 \cup ... \cup E_{n-1}\subseteq E$ where $E_i$ is the set of edges from $v_i$ to $V_i$.
(See Figure~\ref{fig:example_G_sigma} for an example.)  Note that the set $V_i$ -- that is eventually becoming $N^+(v_i)$  and this is the notation we will use for it  --  depends only on the undirected graph $G$, and the  {\em $i$-prefix}  $\sigma_i:=(v_1,\ldots,v_i)$ of $\sigma$. We will also say that $\sigma_i$ is an $i$-prefix in $V$.  If  $N^+(v_{i+1})=\emptyset$, $v_j$ for $j> i$ does not bring in any more edges to $E^\sigma$, so the {\em $i$-suffix} $(v_{i+1},\ldots,v_n)$ of 
$\sigma$ may then be deleted or remain undefined.   We saw that $N^+(v_i)=\emptyset$ if $i>D$, so $\sigma$ can supposed to be an ordered $D$-tuple.

\begin{corollary}\label{cor:co_intervals} Let $G=(V,E)$ be a graph. Then for any ordering $\sigma$ of $V$, $G^\sigma$ is an interval-order subgraph of $G$.  Conversely, any inclusion-wise maximal interval-order subgraph of $G$ is  $G^\sigma$ for some ordering $\sigma$ of $V$. 
\end{corollary}
\smallskip\noindent {\em Proof.}
  Denote by $\overrightarrow{G^\sigma}= (V,\overrightarrow{E^\sigma})$ the digraph obtained by orienting each edge of $G^\sigma$ from its endpoint of smaller index to the one with larger index. Clearly $N^+ (v_i)=V_i$ and by construction, (\ref{eq:chain_property}) is satisfied. Hence, Lemma~\ref{lem:co_intervals} immediately implies the first part of the corollary.
    The converse follows by considering the ordering $\sigma$ satisfying (\ref{eq:chain_property}) that exists by the reverse implication of Lemma~\ref{lem:co_intervals}. 
\qed \smallskip

\section{Interval-orders in line graphs}\label{sec:interval_orders}

In this section, we study the inclusion-wise maximal interval-order subgraphs of the line graph of $K_n$, we can actually list them all (Lemma \ref{lem:maximal_co_interval})!

\medskip

The interval-completion problem for a graph is equivalent,  by complementation,  to finding   interval-order subgraphs in  the complementary graph. The following lemma makes it easier to encounter such subgraphs.  

\begin{lemma}\label{lem:det}  
Let $G=(V,E)$ be a graph, and $\sigma_i = (v_1,v_2, \ldots, v_i)$ an $i$-prefix in  $V$.   If an inclusion-wise maximal interval-order subgraph $H$ belongs to an ordering with prefix $\sigma_i$, then 
\begin{itemize}
\item[(i)] For $u, v\in V \setminus \{v_1,v_2, \ldots, v_i\}$,  $N(u)$ \ $\supseteq N(v)\cap N^+(v_i)$, 
there exists an ordering  $\sigma$ with prefix  $\sigma_i$ such that   $G^\sigma=H$ and $u$ precedes $v$ in $\sigma$. 
\item[(ii)] In every  ordering $\sigma$ with  prefix $\sigma_i$ and $H=G^\sigma$, $\sigma_i$ is immediately followed by all $N_i := \{v \in V \setminus \{s_1, s_2,\ldots, s_i\} :  N(v)\supseteq N^+(v_i)\}$ in arbitrary order.
\item[(iii)]
If $|N^+(v_i)|=2$, then $H$ is one of two possible interval-order graphs. 
\end{itemize}
\end{lemma}

Note that $N_i$ is exactly the set of elements electable to be added to $\sigma_i$ as $v_{i+1}$ so that $N^+(v_{i+1})= N^+(v_{i})$.



\smallskip\noindent {\em Proof of Lemma \ref{lem:det}.}  To show (i) just note: if $v$ precedes $u$, then interchanging $u$ and $v$, the condition of (i) makes sure that neither $N^+(u)$ nor $N^+(v)$ decrease. 

For checking (ii) let  $\sigma$ be an ordering of $V$ with prefix $\sigma_i$ and such that $G^\sigma$ is a maximal interval-order subgraph of $G$.  If $\sigma_i$ is not immediately followed by the vertices of $N_i$, modify $\sigma$ by moving all the vertices of $N_i$ immediately after $v_i$, in arbitrary order.
Clearly, each $N^+(v_i)$ is replaced by a superset, and at least one of them by a proper one: $E^\sigma$ increases, contradicting maximality.
	
	To prove (iii) suppose  $N^+(v_i)=\{x,y\}$, and  $N_x$, $N_y$, $N_{x,y}$ be the vertices in $V\setminus \{v_1,v_2, ..., v_i\}$ adjacent only to $x$, only to $y$,  or to both respectively. Observe that $N_{x,y} = N_i$, so,
	by (ii), $\sigma_i$ is immediately followed by $N_{x,y}$, and then by (i), either an element of $N_x$ or an element of $N_y$  follows, unless both are empty.  
	Applying then (ii) again, the entire $N_x$ or the entire $N_y$ must follow, finishing the proof of the Lemma.
\qed 

\smallskip
We switch now to a higher gear, introducing the key tool to compute the boxicity of $\KG(n,2)$ for $n \geq 5$ (Theorem \ref{thm:box_kneser}), then solve IGC (Lemma~\ref{lem:all_minimal_completion} and Theorem~\ref{thm:minimum_completion}), and finally  compute  the boxicity of complements of line graphs (Theorem~\ref{thm:box_co_line}).  

A solution to these problems necessitates, already for complements of line graphs of small order like the Petersen graph, refined knowledge about their interval completions, that is, about interval-order subgraphs of line graphs. The following lemma establishes that $L(K_n)$ has only  four essentially different maximal interval-order subgraphs.

\begin{lemma}\label{lem:maximal_co_interval} 
For any inclusion-wise maximal interval-order subgraph $H$ of $L(K_n)$ ($n\geq5$), there exist distinct vertices $a, b, c, d, e \in V(K_n)$, uniquely determining the edge-sets $E_{a,b,c,d,e}, E_{a,b,c,d}, F_{a,b,c,d}, F'_{a,b,c,d}\subseteq E(L(K_n))$ such that $E(H)$ is equal to one of these, where
\[|E_{a,b,c,d,e}|= \frac{(n+2)(n-1)}2,\, |E_{a,b,c,d}|= 4(n-1),\, |F_{a,b,c,d}|= |F'_{a,b,c,d}|= 5(n-2).\]
\end{lemma}

The proof consists in presenting the four orderings of $E(K_n)$ determining these sets  as developed in Section~\ref{subsec:ordering}.  After the proof, these sets can be explicitly given with exact formulas. For these formulas, and already for the proof, we denote by $(e,f)$  the edge of $L(K_n)$ between two incident edges $e,f\in E(K_n)$ and, for simplicity, borrow notation from $K_n$ for some edge-sets in $L(K_n)$: 
 \begin{itemize}
\item[-]
for a vertex $v\in V(K_n)$, $Q_v$  denotes the set of $\frac{(n-1)(n-2)}2$ edges of the clique formed in $L(K_n)$ by the $n-1$ edges of $\delta(v) \subseteq E(K_n)$ as vertices of $L(K_n)$;
\item[-]
for an edge $e=uv\in E(K_n)$, $\delta_{uv}$ denotes the set of edges $ef\in E(L(K_n))$, where $f\in V(L(K_n))$ is incident to $u$ or $v$ in $K_n$, that is, $\delta_{uv}$ is the star of center $e=uv\in V(L(K_n))$ in $L(K_n)$, $|\delta_{uv}|= 2(n-2)$;
the set $\delta_{uv^-}$ consists only of the edges $ef$, where $f$ is incident to $u$ in $K_n$, $|\delta_{uv^-}|= n-2$;
\item[-]
for a set $U\subseteq V(K_n)$, $K_U$ denotes the edge-set of $L(K_n[U])$; $K_{u,v,w^-}:=\{(uv,uw), (uv,vw)\}$.  
\end{itemize}

\begin{figure}[ht]
    \centering
    \includegraphics[scale=0.75]{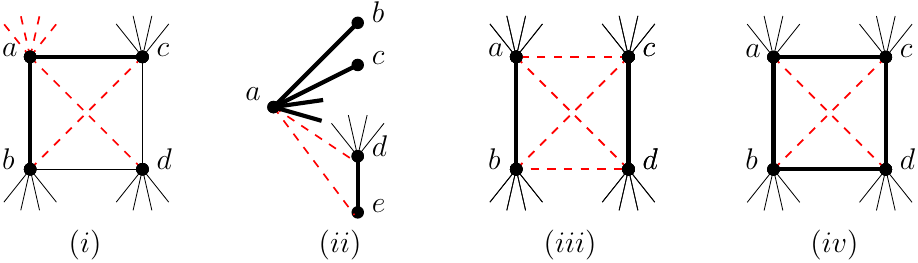}
    \caption{The edges of $V(K_n)$ in $\sigma$;  the edges with fat lines  correspond to $v_1$, $v_2$ in a.left,  and in b.;  to $v_1,\dots, v_{n-2}$ in a.right; and to $v_1, \ldots v_4$ in c.;  the edges with dashed red lines correspond to the vertices in $N^+(v_i)$ $(i=2, n-2, 4)$ respectively. 
    }
    \label{fig:max_int_ord}
\end{figure}

\smallskip\noindent {\em Proof of Lemma~\ref{lem:maximal_co_interval}.}
We define the orderings of $E(K_n)$ giving rise to the four claimed graphs. Let $a,b,c,d,e$ be five different vertices of $K_n$:

First, starting with $v_1=ab$ and $v_2=ac$  adds already $\delta_{ab}\cup \delta_{ac^-} \cup K_{\{a,b,c\}}$ to $E^\sigma$, consisting of $3(n-2)$ edges: $2(n-2)$ edges in $\delta_{ab} \subseteq Q_a \cup Q_b$, $n-3$ edges in $\delta_{ac-} \subseteq Q_a$ (the edge between $ab$ and $ac$ has already been counted) and in addition, the edge between $ac$   and   $bc$, which is  neither in $Q_a$, nor in $Q_b$,  see Figure~\ref{fig:max_int_ord}~$(i)$. Then $v_3,\ldots, v_{n-3}$ is the list of all other edges of $K_n$ incident to  $a$,  except $ad, ae$,  and $v_{n-2}=de$;  we continue by adding the  remaining $n-3$ edges incident to $d$ different from $ad$, in arbitrary order; and finally, $v_{2(n-2)}:= ae$. We added to $E^\sigma$ in this way all edges of $Q_a$, and besides that all the $n-2+1$ edges of $\delta_{ad}\cup K_{\{a,d,e\}}$, symmetrically to the effect of the starting two edges, see Figure~\ref{fig:max_int_ord}~$(ii)$.
Clearly, \[|E^\sigma| =\frac{(n-1)(n-2)}2 + 2(n-1) = \frac{(n+2)(n-1)}2.\] 
The set $E^\sigma$ that we get in this way depends only on the ordered set of the chosen five vertices, so we can denote it by $E_{a,b,c,d,e}$.

Second, let  $v_1=ab$, $v_2=cd$, $v_3=ac$, $v_4=bd$, and continue to define an ordering $\sigma$ with the remaining edges of $\delta(a)$ and of $\delta(d)$ as vertices of $L(K_n)$ in arbitrary order, and then $E^\sigma$. See Figure~\ref{fig:max_int_ord}~$(iii)$ and $(iv)$. Besides the $12$ edges of $K_{a,b,c,d}$ each of $\delta_{ab}$ and $\delta_{ad}$ contain $2(n-4)$  edges of $L(K_n)$, so  $E^\sigma$ has now   $12 + 2(n-4) + 2(n-4)  = 4(n-1)$ edges 
The set $E^\sigma$ depends only on the $4$-tuple $a,b,c,d$ in this order, let us denote it by $E_{a,b,c,d,e}$.

A third kind of interval-order graph arises from  the ordering $\sigma$ with the $4$-prefix defined by $v_1=ab$, $v_2=ac$, $v_3=bd$, $v_4=cd$, leading to $N^+(v_4)=\{ad,bc\}$.  
See Figure~\ref{fig:max_int_ord}~$(i)$ and $(iv)$. Then the ordering continues  either with  the remaining  $n-4$ edges of $\delta(a)$ followed by those of $\delta(d)$, or  the same for $\delta(b)$ followed by $\delta(c)$. 
We denote the corresponding $E^\sigma$ by $F_{a,b,c,d}$, $F'_{a,b,c,d}$ respectively. 

We can now calculate $|F_{\{a,b,c,d\}}|$, $|F'_{\{a,b,c,d\}}|$  for instance by counting the size of the out-neighborhood of each $v_i$:
$|N^+(v_1)| + |N^+(v_2)|= 3(n-2)$, as in the beginning of the proof;
$|N^+(v_3)|=|N^+(v_4)|=2$;  then we add $2(n-4)$ edges of $L(K_n)$ with out-neighborhoods of size $1$, after which the out-neighborhoods are empty, so   $|F_{\{a,b,c,d\}}|= |F'_{\{a,b,c,d\}}|= 3(n-2)+4+2(n-4)=5(n-2)$.

It remains  to  check that any ordering of the vertices defines a subgraph of one of the listed graphs. Due to space constraints, this part of the proof has been postponed to Appendix~\ref{app:proof_lem_4}.
\qed \smallskip

The above proof  shows how $E^\sigma$ can be   explicitly expressed as a function of the vertices  $a,b,c,d,e$ inducing the first edges of $\sigma$ (Figure~\ref{fig:max_int_ord}).  The three formulas are useful to have at hand:

\begin{equation}\label{eq:typeA} \tag{a}
    E_{a,b,c,d,e} = Q_a\cup \delta_{ab}\cup \delta_{ad}\cup K_{\{a,b,c\}}\cup K_{\{a,d,e\}},
\end{equation}
\begin{equation}\label{eq:typeB} \tag{b}
    E_{a,b,c,d}=\delta_{ab}\cup \delta_{ad}\cup K_{\{a,b,c,d\}},
\end{equation}
\begin{equation}\label{eq:typeC} \tag{c}
\begin{aligned}
    F_{a,b,c,d} = \delta_{ab} \cup & \ \delta_{ad}\cup  \delta_{ac^-}\cup K_{\{a,b,c\}}\cup K_{\{a,b,d\}}\cup K_{a,d,c^-}\cup K_{b,c,d^-}. \\
        & \textrm{or $F'_{a,b,c,d}$, where } \delta_{ad} \textrm{ is replaced by } \delta_{bc}.
\end{aligned}
\end{equation}


\section{Boxicity of the Petersen graph and complements of line graphs}\label{sec:petersen}

In this section, we exploit  the description of the inclusion-wise maximal interval-order subgraphs of $L(K_n)$ (Section \ref{sec:interval_orders}),   for proving our main results. We  say that a subgraph  is \emph{ of type (a), (b), or (c) } if its edge-set corresponds to that 
of (\ref{eq:typeA}), (\ref{eq:typeB}), or (\ref{eq:typeC}), respectively; $F_{a,b,c,d}$ and $F'_{a,b,c,d}$ are both considered as of type (c). 

In Section \ref{subsec:proof_main_thm} we establish  Theorem \ref{thm:box_kneser},  that is $\boxi(\overline{L(K_n)}) = n - 2$ for every $n \geq 5$. First, we show the easier upper  bound (Lemma \ref{lem:kneser_upper_bound}). The lower bound, that is, the tightness of the proven upper bound is then proved for $n=5,6$  separately (Lemmas \ref{lem:box_petersen} and \ref{lem:boxlk6}).
These three proofs introduce already the general ideas, but with the difference that for $n\le 6$ the interval-order subgraphs of type  (b) and (c) are larger or equal to those of type (a), so they do play a more essential role.  Then we take on the challenges of the proof for  $n\ge 7$ (Lemma \ref{lem:boxlk7_plus}).
 Section \ref{subsec:box_co_line} makes one more small step for generalizing the results to arbitrary complements of line graphs (without the assumption $G = K_n$), finishing the proofs of Theorem~\ref{thm:minimum_completion} and \ref{thm:box_co_line}. Further possibilities for applying the arguments are discussed in  Section~\ref{sec:conclusion}.


\medskip

\subsection{ Proof of Theorem \ref{thm:box_kneser} }\label{subsec:proof_main_thm}

Theorem 1.1 in~\cite{2023_Caoduro} states that the boxicity of the Kneser-graph $\KG(n,k)$ is at most $n-2$ for $n\geq 2k+1$.
We present a simpler proof for $k=2$ in Appendix~\ref{app:proof_lem_6}.

\begin{lemma}\label{lem:kneser_upper_bound}
	$\boxi(\overline{L(K_n)}) \leq n-2$  if $n\ge 5$.
\end{lemma}
Although the growth of the number of edges is quadratic in  (a), while  it is only linear in (b) and (c), the three sizes are comparable for small values of $n$.
For this reason, the arguments for proving the lower bound in Theorem \ref{thm:box_kneser} for $n\leq 6$ and for $n\geq 7$ are slightly different.

\begin{fact}\label{fact:intersecting_deltas}
	Let $\{a,b,c\}$ and $\{a',b',c'\}$ be two sets of three distinct vertices of $V(K_n)$. 
    Then $\big ( \delta_{ab} \cup \delta_{ac} \big ) \cap \big ( \delta_{a'b'} \cup \delta_{a'c'} \big) = \emptyset$ if and only if $\{a,b,c\}\cap\{a',b',c'\}=\emptyset$.\qed
\end{fact}

\begin{lemma}\label{lem:box_petersen}
	$\boxi(\overline{L(K_5)}) \geq 3$.
\end{lemma}
\smallskip\noindent {\em Proof.}
	Assume for a contradiction that $\boxi(\overline{L(K_5)}) \leq 2$. Let $\{E_1, E_2\}$ be  a $2$-interval-order-cover of $L(K_5)$, and assume that $(V,E_1)$ and $(V,E_2)$ are maximal interval-order subgraphs of $L(K_5)$.
	By Lemma \ref{lem:maximal_co_interval},  $|E_i|=14$  (a), or $|E_i|=15$ (c), or $|E_i|=16$ (b), $(i=1, 2)$. In each of these cases, there are three distinct vertices $a_i,b_i,c_i \in V(K_5)$ such that $\delta_{a_ib_i} \cup \delta_{a_ic_i} \subseteq E_i$ ($i \in \{1,2\}$), (see (\ref{eq:typeA}), (\ref{eq:typeB}), and (\ref{eq:typeC})).
    The set $V(K_5)$ has only $5$ vertices, so $\{a_1,b_1,c_1\}$ and $\{a_2,b_2,c_2\}$ have to intersect, and then $|E_1 \cap E_2| \geq 1$ by Fact~\ref{fact:intersecting_deltas}.
	
	Since $|E(L(K_5))| = 5\binom{4}{2} = 30$, $\{E_1, E_2\}$ forms an interval-order-cover only if $E_1$ and $E_2$ have both at least $15$ edges, and at least one of them has $16$ edges. Assume that $|E_1|=16$ ($E_1$ is of type (b)), and $E_2$ has either $16$ or $15$ edges (it is of type (b) or (c)). 
	In both cases, there are two sets of four distinct vertices in $V(K_5)$ defining the edge-set of $E_1$ and $E_2$.
	These two sets have at least three common vertices, say $\{a,b,c\}$, and $K_{a,b,c}\subseteq E_1$ follows, and also $K_{a,b,c}\subseteq E_2$ if it is of type (b), or $|K_{a,b,c}\cap E_2|\ge 2$  if it is of type (c), (see (\ref{eq:typeB}), and (\ref{eq:typeC})). 
    Either way, $|E_1| + |E_2| - | E_1 \cap E_2|\le 16+15-2 < 30$, contradicting the assumption that  $E_1\cup E_2=E(L(K_5))$.
\qed \smallskip

\begin{lemma}\label{lem:boxlk6}
	$\boxi(\overline{L(K_6)}) \geq 4$.
\end{lemma}
\smallskip\noindent {\em Proof.}
	Assume for a contradiction that $\boxi(\overline{L(K_6)}) \leq 3$. Let $\{E_1,E_2,E_3\}$ be a $3$-interval-order-cover of $L(K_6)$, and assume that $(V,E_1)$, $(V,E_2)$, and $(V,E_3)$ are maximal interval-order subgraphs of $L(K_5)$.
	
	By Lemma \ref{lem:maximal_co_interval}, $|E_i| = 20$  (no matter if it is of type (a), (b), or (c))  and  in all the three cases there are three distinct vertices $a_i,b_i,c_i \in V(K_6)$ such that $\delta_{a_ib_i} \cup \delta_{a_ic_i} \subseteq E_i$.
	Since $|E(L(K_6))| = 6\binom{5}{2} = 60$, $E_1,E_2$ and $E_3$ are pairwise disjoint. Applying Fact~\ref{fact:intersecting_deltas} three times, we deduce that 
    the nine vertices of $a_i,b_i,c_i$ for  $i \in [3] \}$  are all distinct, contradicting $|V(K_6)| = 6$. 
\qed \smallskip

We now finish the proof of  Theorem \ref{thm:box_kneser} with the lower bound for $n\geq 7$.
\begin{lemma}\label{lem:boxlk7_plus}
	$\boxi(\overline{L(K_n)}) \geq n-2$, for any $n\geq7$.
\end{lemma}
\smallskip\noindent {\em Proof.}
	Assume for a contradiction that $\boxi(\overline{L(K_n)}) \leq n-3$. Let $\{E_i : i \in [n-3] \}$ be an $(n-3)$-interval-order-cover of $L(K_n)$, and assume that $\{(V,E_i) : i \in [n-3]\}$ are maximal interval-order subgraphs of $L(K_n)$.
    For each $i \in [n-3]$, there are distinct $a_i, b_i,c_i \in V(K_n)$ such that: 
    \begin{itemize}
        \item [-] $E_i \supset Q_{a_i} \cup \delta_{a_ib_i} \cup \delta_{a_ic_i}$, if $E_i$ is of type (a); or
        \item [-] $E_{i} \supset \delta_{a_{i}b_{i}} \cup \delta_{a_{i}c_{i}}$, if $E_i$ is of type (b) or (c).
    \end{itemize}
    We show that for any possible assignment of $\{(a_i,b_i,c_i) : i \in [n-3]\}$, the number of edges in $\bigcup_{i \in [n-3]} E_i$ is strictly smaller than $|E(L(K_n))|$.
    
	First, observe that \textit{at least $n-4$ interval-order graphs in this cover are of type (a)}.
	Indeed, if there are at most $n-5$ edge-sets of type (a), then, by Lemma~\ref{lem:maximal_co_interval},
    the number of covered edges is at most  $$(n-5) \frac{(n+2)(n-1)}{2} + 10(n-2),$$
	quantity that, for $n\geq 7$, is strictly smaller than $|E(L(K_n))| = n\binom{n-1}{2}$.
	Therefore, the only two cases to consider are: the cover contains $n-3$ edge-sets of type (a), or it contains $n-4$ edge-sets of type (a) and one edge-set of type (b) or (c).

    Then, note that $\sum_{i=1}^{n-3} |E_i|  = |E(L(K_n))|+ \frac{(n-1)(n-6)}{2}$ in the first case, and $\sum_{i=1}^{n-3} |E_i|  = |E(L(K_n))|+ n-6$ in the second case. We prove  $|\bigcup_{i \in [n-3]} E_i|< |E(L(K_n))|$ by showing that, in both cases,  the pairwise intersections of the  $E_i$ $(i\in[n-3])$ sum up to more than $\frac{(n-1)(n-6)}{2}$, and $n-6$, respectively. We omit the details in this limited version.
    
    \qed

\subsection{Interval-completion and boxicity of complements of line graphs}\label{subsec:box_co_line}

Lemma \ref{lem:maximal_co_interval} presents all  the maximal interval-order subgraphs of $L(K_n)$. 
Now, we show how to use this information to generate all the maximal interval-order subgraphs of $L(G)$ for any  graph $G=(V, E)$. Since complementing the maximal interval-order subgraphs of $L(G)$, we get the minimal interval-completions of $\overline{L(G)}$, this will prove Lemma \ref{lem:all_minimal_completion}.   We then derive Theorem~\ref{thm:minimum_completion} and \ref{thm:box_co_line} as   easy applications of this lemma.

\smallskip\noindent {\em Proof of Lemma \ref{lem:all_minimal_completion}.}
    Complete $G=(V, E)$ by adding all of its \emph{non-edges}, to get the complete graph $K_V$. Then  $L(G) = L(K_V)[E]$  is the  subgraph of $L(K_V)$ induced by the vertex-set $E$ of $L(G)$. 
   
    Let $\Ascr$  be the family of edge-sets of all  the at most  $\Oscr (n^5)$ edge-maximal interval-order subgraphs of $L(K_V)$  listed by Lemma~\ref{lem:maximal_co_interval}, and, for each $A\in\Ascr$, denote by $G_A$ the graph $(E(K_V), A)$.  Since interval-order graphs are closed under taking induced subgraphs, the graphs in $\mathcal{B} := \{G_A[E]: A\in\Ascr\}$ (induced by the vertices of $L(G)$ corresponding to the edges of $G$)  are also interval-order graphs.
    
    Each interval-order subgraph of $L(G)$ can be completed to an inclusion-wise maximal interval-order subgraph of $L(K_V)$,  
    so  the inclusion-wise maximal ones among the edge-sets of the graphs in $\Bscr$ are exactly the edge-sets of the maximal interval-order subgraphs of $L(G)$. The edge-set of each can be computed in $\Oscr(n^2)$ time by following the procedure of Lemma~\ref{lem:maximal_co_interval} and by working directly in $K_V$.
\qed \smallskip

Given Lemma \ref{lem:all_minimal_completion} and the notion of $k$-interval-order-covers, the proofs of Theorem~\ref{thm:minimum_completion} and \ref{thm:box_co_line}  follow easily: 

\smallskip\noindent {\em Proof of Theorem \ref{thm:minimum_completion}.}
    Let $G$ be a graph on $n$ vertices.  By Lemma~\ref{lem:all_minimal_completion} all inclusion-wise  minimal interval-completions of $\overline{L(G)}$ can be listed in $\mathcal{O}(n^7)$ time, and we can take the one with a minimum number of edges among them. 
\qed \smallskip

\smallskip\noindent {\em Proof of Theorem \ref{thm:box_co_line}.}
We can assume $k \leq n-3$. 
By Lemma \ref{lem:all_minimal_completion},  the family $\Bscr$ (defined  in the proof of Lemma \ref{lem:all_minimal_completion}) can be computed in   $\Oscr(n^7)$ time and has cardinality $\Oscr(n^5)$. To decide if $\boxi(\overline{L(G)}) \leq k$, one can simply generate all subsets of $k$ distinct elements from $\Bscr$ and return \emph{True} if at least one of these defines a $k$-interval-cover of $L(G)$, and \emph{False}, otherwise.
This algorithm has a running time bounded by $\Oscr(n^7) + \Oscr(n^{5k})f(n,k)$ where $f(n,k)$ is the time spent to check if a $k$-set of $\Bscr$ is an interval-order-cover of $L(G)$.
Since $L(G)$ has at most $n^3$ edges and $k<n$,  $f(n,k) = \Oscr(n^3)$ using an appropriate data-structure. This concludes the proof of the Theorem.
\qed \smallskip

\section{Conclusion}\label{sec:conclusion}

We have explored in this paper the interval-order subgraphs of line graphs and showed that their number can be bounded by a  polynomial of the number of vertices,
and the polynomial solvability of related optimization problems follows. 

 
Some other  connections  extend our arguments concerning the boxicity to other graphs than complements of line graphs, and the new observations have led us to establish the boxicity of some other relevant graphs, about which we provide a short account.   

First, of course, the question of determining the boxicity of line graphs comes up.  Our meta-method dictates to study  the interval-order subgraphs  of the complements of line graphs. The inclusion-wise maximal ones among these are in one-to-one correspondence  with linear orderings (permutations) of the vertex-set. More concretely, for an arbitrary graph $G$ on $n$ vertices,  the boxicity of  $L(G)$  turns out to be  equivalent to a beautiful extremal problem about the permutations of $V(G)$:

We say that a linear ordering $\prec$   of $V(G)$ {\em covers} the pair of vertex-disjoint edges $\{ab, cd\}$  of $G$  if $\max\{a,b\} \prec \min \{c,d\}$,  or $\max\{c,d\} \prec \min \{a,b\}$, where the max and the min concern the linear ordering $\prec$.  We proved that $\boxi (L(G) )= \mathcal{O}(\log n)$ even if $G$ is the largest possible among line graphs, that is,  $G=K_n$. In this case, we also proved  $\boxi (L(K_n) )=\Omega (\log \log n)$.
Further investigations, including algorithms and complexity, are under investigation.  The only general results about the boxicity of these graphs we found in the literature are in \cite{2011_Chandran}.

It is also  tempting to apply   the   meta-method  -- of coloring the edges of the graph so that  each color forms  an interval-order subgraph, and an edge may get several colors -- to some popular  ``Mycielsky graphs'' and  ``Ramsey graphs'' and follow a program of generalizations similar to what we did in this article. An account of these results can be found in \cite[Chapter 4]{Caoduro_thesis}.


We conclude with a few open problems. The first one is from \cite{2023_Caoduro}.

\begin{problem} \cite{2023_Caoduro}
Determine the boxicity of Kneser-graphs $\KG(n,k)$ with $k\geq 3$. 
\end{problem}

Already establishing $\boxi(\KG(7,3))$ is open. The computation of this value could either support or directly refute the intriguing conjecture $\boxi(\KG(n,k)) = n - k$, for any $n\geq 2k + 1$.

Deciding whether the boxicity of the complement of a line graph is at most $k$ can be solved in polynomial time when $k$ is fixed in advance (Theorem \ref{thm:box_co_line}). We did not manage to show that the dependency on $k$ is really necessary.

\begin{problem}\label{prob:co_line_NP_hard}
    Is computing the boxicity of \emph{complements of line graphs}  $\mathcal{NP}$-hard?
\end{problem}

Assuming a positive answer to Problem \ref{prob:co_line_NP_hard}, we ask if the dependency on $k$ can be essentially improved. For instance, \emph{is it FPT?} Problem \ref{prob:co_line_NP_hard} arises for line graphs as well, where even the membership to the class XP is open:

\begin{problem}\label{prob:line_graph}
Can it be decided in polynomial time whether the boxicity of a  \emph{line graph} is $2$, or is it $\mathcal{NP}$-complete?
\end{problem}

It is well-known, and easy to check, that $\boxi(L(K_4)) = 3$. Hence, Problem \ref{prob:line_graph} is interesting only for line graphs of graphs with clique number at most three.

\subsubsection*{Acknowledgment.}~We thank an anonymous referee for numerous helpful suggestions that significantly improved the structure and quality of our presentation.

\printbibliography
\addcontentsline{toc}{chapter}{Bibliography}


\newpage

\appendix

\section{Proof of the missing part of Lemma~\ref{lem:maximal_co_interval} }\label{app:proof_lem_4}

Let $\sigma$ be an ordering of $V$ so that $E^\sigma$ is  inclusion-wise maximal, and let us prove that it has one of the forms (\ref{eq:typeA}), (\ref{eq:typeB}) or (\ref{eq:typeC}):


Let  $v_1:= ab\in E(K_n)$,  and $i(\sigma)$  be the first index of $\sigma$ such that $v_i$, as an edge of $K_n$, is not incident to $a$.  If there is no such a $v_i$, when we  define $i(\sigma):=\infty$.

\begin{figure}[ht]
		\centering
		\includegraphics[scale=0.9]{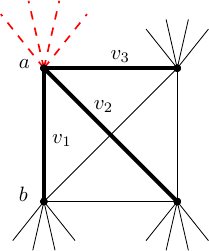}
		\caption{The edges of $V(K_n)$ that can be contained in $\sigma$ in Case 1. The ``legend'' is that of the article, eg. dashed red lines correspond to the vertices in $N^+ (v_3)$.}
		\label{fig:max_int_Case1}
	\end{figure}

	
	\smallskip\noindent {\em Case 1}: $i(\sigma)\ge 4$, that is, $v_1, v_2, v_3\in E(K_n)$ are incident to  $a$ (Fig.~\ref{fig:max_int_Case1}).
 
Denote the two endpoints of $v_i$ by $d$ and $e$.
It follows then from Lemma~\ref{lem:det}~(i), that $N^+(v_{i(\sigma)})=\{ad,ae\}$ can be supposed.

 \smallskip \noindent {\bf Claim}:  $i(\sigma)=n-2$. 

	Indeed, since $N^+(v_{i(\sigma)})=\{ad,ae\}$, we have $i(\sigma)\le n-2$.
  For each edge $af$ (if any), different from $v_1,\ldots, v_{i(\sigma-1)}$ and from $ad, ae$, the neighborhood $N(af)$ contains $N(de) \cap N^+(v_{i(\sigma)-1})$. 
    Hence, applying Lemma~\ref{lem:det}~(i), each edge of this type precedes $de$ and $i(\sigma) = n-2$.

\medskip
    
	
   
    Now, we can finish the ordering using Lemma \ref{lem:det} (iii). The condition $N^+(v_{i(\sigma)})=\{ad,ae\}$ implies that there are only two possible ways to continue: either with the remaining edges of $K_n$ incident to $ad$,  or with those incident to $ae$, in arbitrary order. In this way, we are arriving exactly at $E_{a,b,c,d,e}$ or $E_{a,b,c,e,d}$, respectively.

	\smallskip\noindent{\em Case~2}: $i(\sigma)=3$.
	
\begin{figure}[ht]
		\centering
		\includegraphics[scale=0.75]{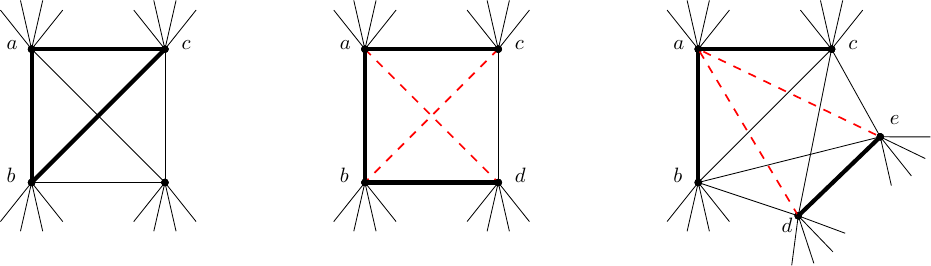}
		\caption{The three edges of $V(K_n)$ that $\sigma$ may start within Case 2. }
		\label{fig:max_int_Case2}
	\end{figure}

	Then $v_2=ac\in E(K_n)$ for some $c\in V(K_n)$ different of $a$ and $b$, 
	and the assumption $i=3$ restricts our choices to the following three possibilities for $v_3$:  
	$v_1$, $v_2$, $v_3$ may form the three edges of a triangle of $K_n$, then $v_3=bc$; or $v_1$, $v_2$, $v_3$ may form a path, that is, $v_3=bd$, or  $v_3=cd$, so that $a,b,c,d$ are four different vertices of $K_n$; or else, $v_3=ed$, where $a,b,c,d,e$ are five different vertices of $K_n$ (Figure~\ref{fig:max_int_Case2}). 
	
	\smallskip {\em Case 2.1}: $v_1$, $v_2$, $v_3$ form the three edges of a triangle of $K_n$ (Figure \ref{fig:max_int_Case2} (left)).
	
	Then $N^+(v_3) = \emptyset$ 
 so any completion of the ordering starting with $v_1, v_2$ will be contained in $\delta_{ab}\cup\delta_{ac}$, a subset of (\ref{eq:typeA}), (\ref{eq:typeB}) or (\ref{eq:typeC}), contradicting maximality.
	
	\smallskip {\em Case 2.2}:  $v_1$, $v_2$, $v_3$ form a path (Figure \ref{fig:max_int_Case2} (center)).
	
	Then $v_3=bd$ (or, equivalently $v_3 = cd$), where $d\in V(K_n)$ is different from $a,b,c$.
	The set of common neighbors of $v_1$, $v_2$, $v_3$ is $\{bc, ad\}$, and by Lemma~\ref{lem:det}~(ii), $v_4=cd$ (or, equivalently $v_4 = bd$, if $v_3 = cd$), and we can either continue the construction of an ordering $\sigma$  with all neighbors in $L(K_n)$ of $ad$,  or all neighbors of $bc$ Lemma~\ref{lem:det} (iii)) arriving respectively at: 
	\[\delta_{ab}\cup \delta_{ad}\cup  \delta_{ac^-}\cup K_{\{a,b,c\}}\cup K_{\{a,b,d\}}\cup K_{a,d,c^-}\cup K_{b,c,d^-},\] 
	\[\delta_{ab}\cup \delta_{bc}\cup  \delta_{ac^-}\cup K_{\{a,b,c\}}\cup K_{\{a,b,d\}}\cup K_{a,d,c^-}\cup K_{b,c,d^-}.\]

    In the first case, we get $E^\sigma=F_{a,b,c,d}$ and in the second case, $E^\sigma=F'_{a,b,c,d}$. Therefore,  Case~2.2 leads anyway to interval-order graphs of the form of (\ref{eq:typeC}).
    
 
	\smallskip {\em Case 2.3}: $v_3=de$, where $a,b,c,d,e$ are different vertices of $K_n$ (Fig. \ref{fig:max_int_Case2} (right)).

	In this case $N^+(v_3)=\{ad, ae\}$, so by Lemma~\ref{lem:det} (i), all the edges incident to $a$, except from $ad$ and $ae$, precede $v_{i(\sigma)}$, that is,  $n-2=i(\sigma)=3$. Hence $n=5$ and, continuing the ordering $\sigma$ as in Case 1, we get the edge-set (\ref{eq:typeA}).

	\smallskip \noindent {\em Case~3}: $i(\sigma)=2$, that is, $\{v_1,v_2\}$ is a matching (Figure \ref{fig:max_int_Case3} (left)).

 A common point of $v_1$ and $v_2$ can be chosen to be $a$ and then $i>2$. So $v_2=cd$, and $a,b,c,d\in V(K_n)$ are different. So $N^+ (v_2)= N_{L(K_n)}(v_1)\cap N_{L(K_n)}(v_2) = \{ad,db,bc, ca\}$, and any choice of $v_3$ corresponding to an edge in $K_n$ incident to $a,b,c$, or $d$ will lead to $|N^+(v_3)|=2$ (Figure~\ref{fig:max_int_Case3}).
	There are two non-equivalent ways of continuing then (similarly to Cases 2.2 and 2.3):
	
	\begin{figure}[ht]
		\centering
		\includegraphics[scale=0.75]{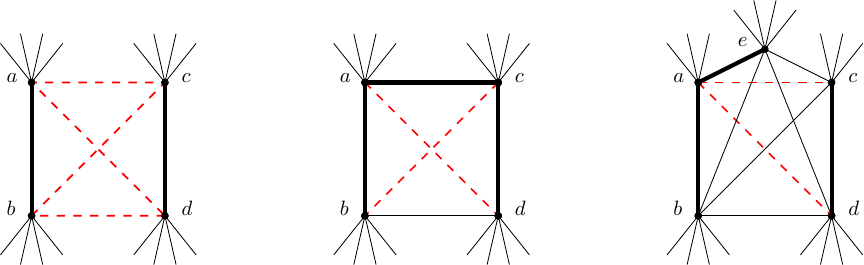}
		\caption{The subcases  of Case~3.}
		\label{fig:max_int_Case3}
	\end{figure}

 {\em Case 3.1}: $v_3=ac$ (Figure \ref{fig:max_int_Case3} (center)).
	
	Then by Lemma~\ref{lem:det} (ii), $v_4=bd$ and the size of the out-neighborhood remains $2$  until adding the next vertex of the ordering, when it is obliged to  decrease  to $1$, and its unique edge can be supposed (by symmetry) to be $ad$. Observing that besides $\delta_{ab}\cup \delta_{ad}$ only the edges of  $K_{\{a,b,c,d\}}$ are covered, we get that  $E^\sigma=\delta_{ab}\cup \delta_{ad}\cup K_{\{a,b,c,d\}}=E_{a,b,c,d}$, so $E^\sigma$ is an edge-set of the form (\ref{eq:typeB}).  
	
	

    \smallskip{\em Case 3.2}: $v_3=ae$,  where $a,b,c,d,e\in V(K_n)$ are different vertices (Figure \ref{fig:max_int_Case3} (right)). 
    
	Then $N^+(v_3)=\{ac,ad\}$, which adds the two edges $(ae,ac), (ae, ad)$ to $E^\sigma$.  By Lemma~\ref{lem:det} (ii), $\sigma$ has to continue with the rest  of the edges incident to $a$, in arbitrary order, finishing with $v_{n-1} = ac$, $N^+(v_{n-1})=\{ad\}$. Applying again Lemma~\ref{lem:det} (ii), we have to complete $\sigma$ adding the remaining $n-3$ edges of $K_n$ that are incident to $d$, again, in arbitrary order. 
 
 We get now the set:

 \[ E^\sigma = \delta_{ab}\cup \delta_{ad}\cup  \delta_{ac^-}\cup K_{\{a,c,d\}}\cup K_{\{a,b,d\}} \cup K_{a,b,c^-}\cup K_{c,d,b^-},\] 
 which is, $F_{a,d,c,b}$, 
concluding again with $E^\sigma$ of the form (\ref{eq:typeC}).
\qed 

\section{Proof of Lemma~\ref{lem:kneser_upper_bound} }\label{app:proof_lem_6}
    In the proof, we denote the vertices of $V(K_n)$ by $v_1, v_2, \ldots, v_n$. 
    The edge-set of $L(K_n)$ is $\bigcup_{i=1}^n Q_{v_i}$. We show that it can be covered using $(n-2)$ interval-order graphs of type (a).

	For $i \in [n-2]$, let $G_i = (V,E_i)$ be an interval-order subgraph of $L(K_n)$ with edge-set $E_i \supseteq Q_{v_i} \cup \delta_{v_iv_{n-1}} \cup \delta_{v_iv_{n}}$, $(i=1,\ldots, n-2)$. The union of these obvious covers
$Q_{v_i} \subset E_i$ for all indices, except possibly for $i=n-1$ or $i=n$.  However, observe that 
    $Q_{v_{n-1}} \subseteq \bigcup_{i \in [n], \  i \neq n-1} \delta_{v_iv_{n-1}}$, and similarly, 
    $Q_{v_n} \subseteq \bigcup_{i \in [n-1]} \delta_{v_iv_{n}}$.
    These sets are also covered since for $i \in [n-2]$, $\delta_{v_iv_{n-1}}, \delta_{v_iv_{n}} \subset E_i$, and $\delta_{v_{n-1}v_{n}}$ consists of the $2(n-2)$ edges $\{(v_iv_{n-1},$ $v_{n-1}v_{n}) : i \in [n-2] \}$ and $\{(v_iv_{n}, v_{n-1}v_{n}) : i \in [n-2]\}$ which are contained in $\delta_{v_iv_{n-1}}$ and $\delta_{v_iv_{n}}$, respectively.
\qed \smallskip
\end{document}